\title{A proof that $\sqrt{s}$ for $s$ not a perfect square is simply normal to base 2 }
\author{Richard Isaac}
\documentclass[11pt]{article}
\newtheorem{theorem}{Theorem}
\newtheorem{corollary}{Corollary}
\newtheorem{lemma}{Lemma}

\date{   }
\begin{document}
\maketitle
\begin{abstract}
\noindent  Let   $\omega$ have dyadic expansion $.x_{1}x_{2} \cdots$. Put  $f_{n}(\omega)=(x_{1}+x_{2}+\cdots+x_{n}) /n$. If $\lim f_{n}(\omega) \rightarrow 1/2$, $\omega$ is called simply normal to base 2.  We give a proof  that numbers of the form $\sqrt{s}$ for $s$ not a perfect square have dyadic expansions which are simply normal.   

Let $n_{k}$ be an arbitrary subsequence and define $f(\omega)=\limsup f_{n_{k}}(\omega)$. We define two expansions, $\omega_{1}$ and $\omega_{2}$ such that  $f(\omega_{1})=a$ and $f(\omega_{2})=1-a$ where $1-a=f(\sqrt{s})$. Let $\nu=\omega^{2}$ and define $h(\nu)=f(\omega)$. 
 It is observed that the $\nu_{i}$, the squares of the $\omega_{i}, i=1,2$, have the same tail, that is, there exists an index $n$ such that the digits at indices $n+j, j \geq 0$ are the same for both expansions. From the above, $h(\nu_{2})=1-h(\nu_{1})$. The proof is completed in section 2.3 by showing $h$ at the $\nu_{i}$ does not depend on any initial segment of its digits. This part depends on basic probability involving the notions of independence and conditional expectation. 
 \footnote{{\it  AMS 2010 subject classifications}, 11K16  \\
{\it \hspace*{1em} Keywords and phrases}. normal and simply normal numbers, tail of an expansion, tail                
\hspace*{1.5in} function, independent random variables, conditional expecta- 
\hspace*{1.5in} tion, see section 2.3.2.. }\\ 
  \end{abstract}
\section{Introduction} 
 \subsection{ The problem and its setting}
A number is {\it simply normal to base $b$} if
its base $b$ expansion has each digit appearing with average
frequency tending to $b^{-1}$. It is {\it normal to base $b$} if its
base $b$ expansion has each block of $n$ digits
appearing with average frequency tending to $b^{-n}$. A number is called
{\it normal}
if it is normal to base $b$ for every base. For a more detailed
introductory discussion we refer to chapter 8 of \cite{ni:ir} or section 9.11 of \cite{hw:tn}.
The most important theorem about normal numbers is the celebrated result
(1909) of E. Borel in which he proved the normality of almost all
numbers with respect to Lebesgue measure (for a proof see section 9.13 of \cite{hw:tn}; an elegant probability proof using the strong law of large numbers appears in \cite{jl:mt}, p. 43). 

Borel's theorem left open the question,
however, of identifying specific numbers as normal, or even
exhibiting a common irrational simply normal number.  \markright{INTRODUCTION}
In this paper we exhibit a class of numbers simply normal to the base 2.
More precisely, we prove
\begin{theorem} \label{snb2} Let $s$ be a natural number which is not a perfect
square. Then the dyadic (base 2) expansion of $\sqrt{s}$ is simply normal.  \end{theorem}
Consider irrational  numbers $\omega$ in the closed unit interval $\Omega$, and represent the dyadic expansion of $\omega$ as
\begin{equation} \label{onee}
\omega=.x_{1}x_{2} \cdots ,\hspace{.4in} x_{i}=0 \mbox{ or } 1.  \label{expom}  \end{equation}
This expansion is unique because of the irrationality of $\omega$.

  \subsection{Tail functions, coordinate averages}
Given any  expansion $ \omega=.x_{1}x_{2}\cdots$ and any natural number $n$, the sequence of digits $x_{n}, x_{n+1}, \cdots$ is called a {\it tail} of the expansion. Two expansions are said to have the same tail if there exists $n$ so large that the tails of the sequences from the $n$th digit are equal  (that is, the digits at  indices $n+j$ for  $j\geq 0$ are the same for both expansions). 
 The function $g$  is called a {\it tail function}  if, whenever $\alpha_{1}$ and $\alpha_{2}$ have the same tail,  $g(\alpha_{1})=g(\alpha_{2})$. This means that for every  natural number $n$ there is a function $g_{n}$ defined on the sequence $x_{n},x_{n+1}, \cdots$ such that 
 \begin{equation} \label{ttail} g(x_{1},x_{2},\cdots)=g_{n}(x_{n},x_{n+1} \cdots). \end{equation} 
 
 The average
\begin{equation} \label{aver}
f_{n}(\omega)= \frac{x_{1}+x_{2}+\cdots+x_{n}}{n} \end{equation}
is the relative frequency of 1's in the first $n$ digits of the expansion of $\omega$. Simple normality for $\omega$ is the assertion that $f_{n}(\omega) \rightarrow 1/2$ as $n \rightarrow \infty$.
 Let $n_{k}$ be any fixed subsequence and define \begin{equation} \label{ff}
f(\omega)= \limsup_{k \rightarrow \infty} f_{n_{k}}(\omega) \end{equation}
where $f$ depends on the subsequence $n_{k}$.
 Note that  \begin{equation} \label{are}
 f(\omega)=  \limsup_{k} \frac{x_{1}+x_{2}+ \cdots x_{n_{k}}}{n_{k}}=\limsup_{k} \frac{x_{r}+x_{r+1}+ \cdots x_{n_{k}}}{n_{k}} \end{equation}
no matter how large fixed $r$ is, so that $f$ is a tail function. Moreover, because of the second equality $f$ is also invariant, that is, $f(T \omega)=f(\omega)$ where $T$ is the shift transformation taking $x_{1},x_{2},\cdots$ into $x_{2},x_{3},\cdots$.
 \subsection{$f_{n}$  and $f$  on $\omega$ written as  functions $h_{n}$ and $h$ on $\nu=\omega^{2}$} \label{ffo}
In the arguments to follow the square $\nu=\omega^{2}$ of an expansion $\omega$ will play an important part. Let us put
\begin{equation} \label{xu}
 \omega^{2}=\nu=.u_{1}u_{2} \cdots, \hspace{.4in} u_{i}=0 \mbox{ or } 1. \end{equation} 
Then $\omega$ and $\nu$ uniquely determine each other. So the average $f_{n}(\omega)$, defined in terms of the $x$ sequence of relation~\ref{expom}, can also be expressed as a function $h_{n}(\nu)$ of the $u$ sequence of relation~\ref{xu}.  This relationship has the simple form  
$f_{n}(\omega)=f_{n}(\sqrt{\nu})=h_{n}(\nu)$. Let $n_{k}$ and $f$ be as defined in section 1.2. Define $h(\nu)= \limsup_{k} h_{n_{k}}(\nu)$; then clearly 
$f(\omega)=h(\nu)$.  

\subsection{Idea of proof}
Here is a guide to the basic line of reasoning in the proof of theorem~\ref{snb2}.  Two expansions $\omega_{1}$ and $\omega_{2}$ are defined such that $f(\omega_{1})=a$ and $f(\omega_{2})=1-a$ where $1-a=f(\sqrt{s})$ (see relation \ref{ff}). Let $\nu_{i}=\omega^{2}_{i}, \,  i=1,2$.   Looking at the expansions of these numbers, we see that the $\nu_{i}$ have the same tail. Moreover, the tails of the $\nu_{i}$ and  the  $\omega_{i}$ are related in a simple functional way so they contain the same information. This allows a calculation of $h$ at the $\nu_{i}$ using tail data of the $\nu_{i}$. Section 2.3 digs a little deeper into the relationship between the tails of the $\nu_{i}$ and  the  $\omega_{i}$. Expressing this relationship with two formulas, relations~\ref{ux1} and \ref{urx2}, we show these formulas contain all information about the expansions and the values of $h$. Finally, we use a probability argument using independence and conditional expectation to show $h$ at the  $\nu_{i}$ is just a function of its tails. 
  
\section{Proof of theorem~\ref{snb2}} 
     \label{pts2} 
\subsection{Definition of $\omega_{i}$ and $\nu_{i}$ for $i=1,2$ and determination of their tails} 
     Let $s$ be a natural number which is not a perfect
square, and let $l$ be any positive integer such that $2^{l} > s$.
Define the numbers
\[\omega_{1}=1-(\sqrt{s}/2^{2l})  \]
and
\[ \omega_{2}=(\sqrt{s}-1)/2^{l}.  \]
 
 The numbers $\omega_{i}$ are less than 1 and their squares are  given by
\begin{eqnarray} \label{twopts}
    \nu_{1}=1+s(2^{-4l})-(2^{-2l+1} \sqrt{s}) \mbox{  and  }
 \nu_{2} = (s+1)2^{-2l}-(2^{-2l+1} \sqrt{s}). \end{eqnarray}
Let us study what the expansions of the numbers in relation~\ref{twopts} look like.
The  expansions of the rational terms $1+s(2^{-4l})$ and $(s+1)2^{-2l}$
 have only a finite number of non-zero digits.
 The  expansion of the term $2^{-2l+1} \sqrt{s}$\, 
is obtained from the expansion of
$\sqrt{s}$ by shifting the ``decimal'' point $2l-1$ places to the
left. To get each of the values in relation~\ref{twopts}, this term
must be subtracted from each of the larger rational terms which
have terminating expansions. Let the non-integer part of $\sqrt{s}$ be $.s_{1}s_{2}\cdots$. 
Then the expansions of the squares  look something like the following example:
\begin{eqnarray} \label{sub}
A,1- .000\cdots s_{1}s_{2}s_{3}\cdots &=&A, 0 11\cdots -.000\cdots s_{1}s_{2}s_{3} \cdots\\
& =& B, (1-s_{i})\, (1-s_{i+1}) \cdots \nonumber \end{eqnarray}
Here $A,1$  represents the rational term, a finite sequence of digits $A$ followed by a terminal 1. Subtracted from this is the expansion of the shifted $\sqrt{s}$.  
To perform the subtraction, change the rational term into an infinite expansion by replacing the 
terminal 1 with 0 and all following  0 digits to 1's. The result is an expansion having an initial segment of  digits B followed by a tail of the sequence  $1-s_{i}$. It is important to note that because the shift of $\sqrt{s}$ is the same for both $\nu_{1}$ and $\nu_{2}$, if, for fixed index $i$, the term  $1-s_{i}$ on the right hand side of relation~\ref{sub} appears at index $n$ of the expansion of $\nu_{1}$, then  $1-s_{i}$ also appears at index $n$ of the expansion of $\nu_{2}$. \\
{\bf Definition} Given two sequences $\alpha=.a_{1}a_{2}\cdots$ and $\beta=.b_{1}b_{2}\cdots$, we say the tail of $\alpha$ has tail a tail of $\beta$ to mean there exist natural numbers $N$ and $M$ such that $a_{N+j}=b_{M+j}, j \geq 0.$ \\
 Conclude from relation~\ref{sub} and the subsequent analysis that: \\
{\it the expansions of  $\nu_{1}$ and  $\nu_{2}$ have the same tail. This tail is a tail of the sequence $1-s_{i}$}. 

Examine now the expansions of the $\omega_{i}$. The number $\omega_{1}$ has the same form as the numbers of relation~\ref{twopts}: a shift of $\sqrt{s}$ subtracted from a rational term. So the expansion of $\omega_{1}$ has  tail a tail of the sequence $1-s_{i}$. The expansion of $\omega_{2}$, however, involves a rational term subtracted from a  shifted $\sqrt{s}$. It is easy to see that the result is an expansion with tail a tail of the sequence $s_{i}$. To summarize: \\
{\it the expansion of $\omega_{1}$ has tail a tail of the sequence $1-s_{i}$ and the tail of $\omega_{2}$ has tail a tail of the sequence $s_{i}$.}
\subsection{Calculation of $f(\omega_{i})$;  calculation of  $h(\nu_{i})$ from the tails of  $\nu_{i}$ and their relation to the tails of $\omega_{i}$}
Let $n_{k}$ be any subsequence such that $ \lim_{k \rightarrow \infty} f_{n_{k}}(\omega_{1})$ converges to a limit $f(\omega_{1})=a$. Since $f$ is a tail function, the tail of  $ \omega_{1}$ shows that $a= 1-f(\sqrt{s})$ and the tail of $ \omega_{2}$ is a tail of $s_{i}$ and so
$ \lim_{k \rightarrow \infty} f_{n_{k}}(\omega_{2})$ converges to $1-a=f(\sqrt{s})$ (see relation~\ref{are} and section 2.1). 
For general $\omega$ define $f(\omega)=\limsup _{k \rightarrow \infty} f_{n_{k}}(\omega)$. 
By section~\ref{ffo}
 \begin{equation} \label{s1}
 a= f(\omega_{1})=h(\nu_{1}) \hspace{2em} \mbox{ and  } \hspace{2em} 1-a =f(\omega_{2})=h(\nu_{2}). \end{equation}
 
   Let $x_{ij}$ be the coordinates of $\omega_{i}$ and $u_{ij}$ the coordinates of $\nu_{i}$ for $i=1,2$. 
 From the tail properties described in section 2.1 we can take $r$ so large such that there are indices $n$ and $k$  with
 \begin{eqnarray} \label{fandh}
 f(x_{11}, x_{12}, \cdots x_{1,n-1},1- s_{r},1- s_{r+1},\cdots)=f(\omega_{1})=h(\nu_{1})= && \\
 h(u_{11}, u_{12}, \cdots u_{1,k-1}, 1-s_{r}, 1-s_{r+1},\cdots). && \nonumber \end{eqnarray}
 Observe that the tails of $\nu_{1}$ are the same as those of $\omega_{1}$ up to a shift.  We  show that  $h(\nu_{1})$ can be calculated using the tail data of $\nu_{1}$. Take a typical tail $1- s_{r},1- s_{r+1},\cdots$ of $\nu_{1}$. Since $f$ is a tail function with respect to the $x$ variables we have
\[ a= \lim_{k \rightarrow \infty} ((1-s_{r})+ (1-s_{r+1})+  \cdots +(1-s_{n_{k}}))/n_{k} =f(\omega_{1})=h(\nu_{1}).  \]
Conclude that  one is able to calculate the value $a=h(\nu_{1})$ using the tail data  for $\nu_{1}$ and its relation to the tail of $\omega_{1}$.

The argument for  the pair $(\omega_{2}, \nu_{2})$ is similar. 
 We may assume $r$ has been chosen so large that index $k$ begins a common tail and there is an index $m$,  so that 
\begin{eqnarray} \label{f2andh2}
 f(x_{21}, x_{22}, \cdots x_{2, m-1}, s_{r}, s_{r+1},\cdots)=f(\omega_{2})=h(\nu_{2}) && \\
 =h(u_{21}, u_{22}, \cdots u_{2,k-1}, 1-s_{r}, 1-s_{r+1},\cdots). && \nonumber \end{eqnarray}
First note  that by section 2.1,  $\nu_{2}$ has the same tail as  $\nu_{1}$. Also, up to a shift the tail digits of $\nu_{2}$ and those of $\omega_{2}$  add to 1. So given a typical tail $1- s_{r},1- s_{r+1},\cdots$ of $\nu_{2}$, the tail $s_{r}, s_{r+1}, \cdots$ is also given,
 and then the tail function property of $f$  implies
\[1- a= \lim_{k \rightarrow \infty} (s_{r}+ s_{r+1}+  \cdots +s_{n_{k}})/n_{k} =f(\omega_{2})=h(\nu_{2}).  \]
 Therefore $h$ can be calculated at $\nu_{2}$ from the tail data and its relation to the tail of $\omega_{2}$. 
 \subsection{Proof that $h(\nu_{1})=h(\nu_{2})$, and conclusion} 
  We would like to assert that at the $\nu_{i}$, $h$ is a function only of the tails of the $\nu_{i}$, that the leading digits (i.e., $u_{ij},\,  j<k$) provide no added information about $h$. Since the tails of the $\nu_{i}$ are equal it would then follow that $a=h(\nu_{1})=h(\nu_{2})=1-a$ and $a=1/2$. 
  
Although $h$ is not a tail function, $h$ is related to $f$ which is a tail function; moreover the tails of $\omega_{i}$ and $\nu_{i}$ are intimately related,  as seen in relations~\ref{fandh} and \ref{f2andh2}. These relations are what allowed us to perform the calculations in section 2.2.

Consider the digits $x_{j}$ and $u_{j}$  of relations \ref{expom}   and \ref{xu} as variables, so that we can write for  $\nu_{1}$  (see relation~\ref{fandh})
  \begin{equation} \label{ux1}
u_{k+j}=u_{1,k+j}=1-s_{r+j}=x_{n+j} \hspace{2em}   j \geq 0  \end{equation}
and for  $\nu_{2}$    (see relation~\ref{f2andh2})
\begin{equation} \label{urx2}
u_{k+j}=u_{2,k+j}=1-s_{r+j}=1-x_{m+j}  \hspace{2em}   j \geq 0. \end{equation} 
From these relations, given any $u$ tail of $\nu_{1}$ or $\nu_{2}$,  one can calculate the associated $x$ tail, and conversely. 
\subsubsection{Two results needed for the proof of theorem~\ref{indep2}}
\noindent {\bf Definition} An initial segment of length $n$ of an expansion $\omega$ is the finite expansion $x_{1} x_{2} \cdots x_{n}$ of the first $n$ digits of $\omega.$    \\

We will need to make some elementary observations about the relationship between the digits in the expansion of $\omega$ and those in the expansion of  $\nu=\omega^{2}$.
\begin{lemma} \label{existN}  Let $\omega$ and $\nu=\omega^{2}$  be irrational expansions. Let $k$ be a fixed positive integer. Then there exists a positive integer $N$ (depending on $\omega$) such that the digits $u_{1}, u_{2}, \cdots u_{k}$ of $\nu$ can each be written as a function of the digits 
$x_{1}, x_{2}, \cdots x_{N}$ of $\omega$. \end{lemma}
Proof: The expansions are unique due to irrationality. Let $\omega_{n}$ be the initial segment of length $n$ of $\omega$.
Then $ \omega_{n} \rightarrow \omega$ in the standard metric topology and then so does $\omega_{n}^{2} \rightarrow \nu$.  It must be shown that there is an integer $N$ such that $\omega_{n}^{2},\, n>N$ always produces the given $k$ digits in its approximation to $\nu$. The convergence  $\omega_{n}^{2} \rightarrow \nu$  implies that the $i$th digit of $ \omega_{n}^{2}$ converges to the $i$th digit of $\nu$, precluding any infinite oscillation. Take $N$ large enough so that the initial segment of length $k$ has converged. {\bf Q.E.D.}\\

 How much information is known about $\omega$ and $\nu$ if one only knows the tails of both? The next result provides an answer.
  \begin{lemma} \label{recon}
For any expansions $\omega$ and $\nu= \omega^{2}$,  $\omega$ and $\nu$ can be reconstructed given their tails. \end{lemma}
Proof: Given a tail of $\omega$, there are only a finite number of ways to fill in its initial segment and at least one of these must have a square with tail $\nu$. Suppose $\omega$ and $\omega^{*}$ have the same tail. Going out far enough in the sequence it may be assumed that 
$\omega=A+B$ and $\omega^{*}= A^{*}+B$ where $A \neq A^{*}$, the $A$ terms rational  and $B$ irrational. The difference in these squares is $A^{2}-(A^{*})^{2}+2(A-A^{*})B = \nu-\nu^{*}$ an irrational number. But if  $ \nu$ and $\nu^{*}$ had the same tail $\nu-\nu^{*}$ would be a finite rational. Therefore distinct $\omega$ expansions  with the same tail must have $\nu$ expansions with different tails.   {\bf Q.E.D.}\\
\begin{theorem}  \label{indep2} The function $h$ at the $\nu_{i}$ is only a function of its tails and the relations~\ref{ux1} and \ref{urx2}; it does not depend on any initial segment of its digits. Consequently 
$h(\nu_{1})=h(\nu_{2})$. \end{theorem}

 To see that $h$ at the $\nu_{i}$ does not depend on initial segments, a  probability argument is used. The basic idea is this: for an appropriate probability distribution on the $\omega$ sample space it will be seen that initial segments of  $\nu_{i}$ are independent of tail sets of  $\nu_{i}$ sufficiently far out in the sequence. By relations~\ref{ux1} and \ref{urx2}, $h$ is a function of tails of $\nu_{i} $ and  $\omega_{i}$. Independence implies that the conditional expectation of $h$ given the initial and tail conditions only depends on the tail conditions, and this leads to the  result. 
 
  For an introduction to the probability used below see, e.g., \cite{ri:pp} (for an elementary treatment without measure theory),  \cite{lb:tp} or \cite{rd:te} (for  more advanced studies).  Before we get to the proof of the theorem we present a list of some of the basic definitions and notations from probability used in our arguments. 
\subsubsection{ A list of some probability basics  }  
Note: This section is for reference. It may be skipped until needed. Although integrals are used in some of the definitions below for generality, applications in our arguments will always be to discrete situations where integral simplifies to addition.
\begin{enumerate}
\item Probability space: a triple $( \Omega, P, \mathcal{S})$, where $\Omega$ is a state space, and $P$ a probability on a permissible class of sets $\mathcal{S}$ (closed under countable unions, intersections and complements).
\item The expectation  of a random variable  $X$, denoted $EX$ (or $E_{P} X$): the average of values of $X$ using the probability $P,\,\, \int x\, dP$.
\item The conditional probability of $X$ given a set $S$: defined as the distribution of $X$ given by 
\[P(X \epsilon A \mid S )= \frac{P([X \epsilon A] \cap S)}{P(S)}.\]  This represents an updating of $P$ to a new probability distribution reflecting the given information $S$. Note that  $P(X \epsilon S \mid S )=1$. If given several sets $S_{1}, S_{2}, etc.$, the notation $P(X \epsilon A \mid S_{1}, S_{2}, etc. )$ means $S$ is replaced by the intersection of $S_{1}, S_{2}, etc.$  in the preceding rule. 
\item Independence: Sets $A$ and $B$ are said to be independent if $P(A \mid B)=P(A)$, equivalently $P( A \cap B)= P(A) \cdot P(B)$. The interpretation is that $B$ adds no new information to update $P(A).$
\item The conditional probability of $X$ given a random variable  \newline  $Y:   P(X \epsilon A \mid Y). $  This defines a family of distributions depending on given values of $Y$. If  given $Y \epsilon B$ or $Y=b$, for example,  a set is defined and one uses the preceding rule. The definition extends to the multidimensional case, e.g. $P(X \epsilon A \mid Y_{1}, Y_{2}, \cdots)$.
 Random variables $X$ and $Y$ are independent if the sets $X \epsilon A$ and $Y \epsilon B$ are independent for all $A, B$. Given random variables  $X_{1}, X_{2}, \cdots$ and $Y_{1}, Y_{2}, \cdots$, if the $X$ variables are independent of the $Y$ variables, then any function of the $X$ variables is independent of any function of the $Y$ variables.

\item The conditional expectation of $X$ given $Y$: this is a function of Y defined by taking the averages of $X$ over the family of conditional probabilities given above. Thus \[E(X \mid Y=b)=\int X\, d P(X \mid Y=b).\]	
The conditional expectation, as a function of $Y$, is also a random variable with value 
$E(X \mid Y)(\omega)=E(X \mid  Y=Y(\omega)\,)$
\item Alternate definition of conditional expectation: $E(X \mid Y) (\omega)$ is a random variable  such that if $B$ is any set defined in terms of $Y$ (B is said to be $Y$ measurable) then \[ \int_{B} E(X \mid Y) (\omega)\, dP =\int_{B} X(\omega)\, dP\]
 That this definition is equivalent to the previous  one is an easy exercise. The conditional expectation is defined up to sets of probability zero, so there can exist different versions.
 \item Let $X_{i}$ be a sequence of random variables on the probability spaces $(\Omega_{i}, P_{i},  \mathcal{S}_{i})$  where $\Omega_{i}$ is a subset of the reals. Then the product of the measures $P_{i}$ on product space relative to the product class of sets is a probability space on which the $X_{i}$ are independent random variables. 
  \end{enumerate}
 \subsubsection{Construction of the probability space }
Proof:  To show $h$ independent of initial segments, the first step is to construct a probability space from two other probability spaces in a certain way. First define 
$(\Omega, \pi_{1}, \mathcal{B})$, where $\Omega$ is the set of expansions in the unit interval and  $\mathcal{B}$  the Borel sets of the reals. Define the probability distribution $\pi_{1}$ assigning positive probability to $\omega_{1}$ as follows. Consider a convergent infinite product $\Pi p_{j} \rightarrow p, 0<p<1$. Assign the mass $P(x_{j}=x_{1j})=p_{j}$. Let the product probability be $\pi_{1}$; under it the $x_{j}$ are now independent random variables,  $\pi_{1}(\omega_{1})=p>0$ and the total mass is concentrated on expansions with tail the same as that of $\omega_{1}$. Clearly the $u$ variables, as functions of the $x$ variables, are also random variables. 

The second probability space $(\Omega, \pi_{2}, \mathcal{B})$ follows the same pattern. This space has the same set of outcomes and Borel sets but differs in the probability distribution. This time construct a distribution  $\pi_{2}$ assigning positive probability to $\omega_{2}$. Let the  convergent infinite product $\Pi q_{j} \rightarrow q, 0<q<1$. Assign the mass $P(x_{j}=x_{2j})=q_{j}$. The  product measure $\pi_{2}$, makes the  $x_{j}$  independent random variables,  $\pi_{2}(\omega_{2})=q>0$ and the total mass is concentrated on expansions with tail the same as that of $\omega_{2}$.

For the final step toss a fair coin, say, whose outcomes are head, $Z=H$ or tail, $Z=T$.  Let  $\Omega_{1}$ be the sample space consisting of  all  pairs $(Z, \omega)$. Assign the probability $P$ by setting $P(Z, \omega)= \frac{1}{2} \pi(\omega)$ where $\pi= \pi_{1}$ or $\pi_{2}$ according to whether the value of $Z= H$ or $T$, respectively. Note that  $x$ and $Z$ variables are not independent: information about $Z$ changes the probability of the $x$'s. The process can be thought of as a kind of game: the player tosses a coin. If heads comes up, the player enters the universe governed by $\pi_{1}$, otherwise entering that of $\pi_{2}$. The point $\omega$ has probability given by
\[P( \omega)=\frac{1}{2} \pi_{1}(\omega) + \frac{1}{2} \pi_{2}(\omega ) \]
 and
\[P(\omega_{1})= \frac{1}{2}\pi_{1}(\omega_{1})=\frac{p}{2}>0 \mbox{  and  } P(\omega_{2})= \frac{1}{2}\pi_{2}(\omega_{2})=\frac{q}{2}>0. \]

\subsubsection{ Proof of theorem 2} We already know that $h$ at the $\nu_{i}$ is a function of its tails and relations~\ref{ux1} and \ref{urx2}. To show independence of initial segments first recall
 lemma~\ref{existN}, so that for  $\omega_{1}$ there exists a positive integer $N_{1}$ such that the digits 
$u_{11}, u_{12}, \cdots u_{1k}$ can all be expressed as a function of the digits $x_{11}. x_{12}, \cdots x_{1, N_{1}}$.
For $\omega_{2}$ let $N_{2}$ be the corresponding integer. If we let $N=\max(N_{1},N_{2})$ then the initial segments of length $k$ of both $\nu_{1}$ and $\nu_{2}$ are expressible as functions of the variables $x_{1}, x_{2}, \cdots x_{N}$.   
From relations~\ref{ux1} and \ref{urx2} an index $k=t$ can be chosen so large that the indices $n$ and $m$ for $x$ in those relations satisfy $\min(n,m)=M>N$. In this case the variables in the tail $ u_{t}. u_{t+1}, \cdots $ are functions of the variables $x_{M}, x_{M+1}, \cdots$ for both $\nu_{1}$ and $\nu_{2}$. 

Let the term $1-s_{d}$ correspond to $u_{t}$ in the representations of relations~\ref{ux1} and \ref{urx2}. Then the tail of both $\nu_{1}$ and $\nu_{2}$ starting from index $t$ can be written \begin{equation}
B_{t}=\{u_{t}=1-s_{d},\, u_{t+1}= 1-s_{d+1}, \cdots\}.  \end{equation}
Let  $A_{i,k}=\{ u_{1}=u_{i1}, u_{2}=u_{i2}, \cdots  u_{k}=u_{ik}\},\,  i=1,2. $ 
\begin{lemma}  \label{abp}  \begin{enumerate}
\item The sets $A_{i,k}$ are independent of the set $B_{t}$  with respect to both $\pi_{1}$ and $\pi_{2}$. 
\item  The sets $A_{i,k}$ are independent of the sets $[h(\nu)=h(\nu_{j})] \cap B_{t} \},\, i=1,2;\, j=1,2$ with respect to both $\pi_{1}$ and $\pi_{2}$. \end{enumerate} \end{lemma}
Proof: From the above it is seen that $A_{i,k}$ is defined in terms of  $x$ variables of index at most $N$ and $B_{t}$ defined on $x$ variables of index larger than $N$. The $x$ variables are independent with respect to both $\pi_{1}$ and $\pi_{2}$, proving the first assertion. For the second assertion note that the set $\{h(\nu)=h(\nu_{j}\}$ is defined in terms of 
  $x$ variables of index larger than $N$ by relations~\ref{ux1} and \ref{urx2}. Since $B_{t}$ is also defined in terms of  $x$ variables of index larger than $N$, so is their intersection. Therefore this set is independent of $A_{i,k}$.  {\bf Q.E.D.}
  \begin{lemma}   \label{prob0} Let $\nu_{3}$ be an expansion different from both $\nu_{1}$ and $\nu_{2}$ having tail $B_{t}$. Then $P(\nu_{3})=0.$ \end{lemma}
Proof: From lemma~\ref{recon} it follows that no expansion $\omega_{3}= \sqrt{ \nu_{3}}$ with the same tail as either $\omega_{1}$ or $\omega_{2}$ could have $\nu_{3}$ with tail $B_{t}$. But then $P(\nu_{3})=P(\omega_{3})=0$ since $P$ only puts positive mass on points with tail the same as those of $\omega_{1}$ and $\omega_{2}$. {\bf Q.E.D.}\\

The probability 
\[P( h(\nu) \epsilon S \mid u_{1}, u_{2}, \cdots u_{k}; B_{t})\] 
only assigns positive measure to the expansions $\nu_{1}$ and $\nu_{2}$ by lemma~\ref{prob0}. The vector  $u_{1}, u_{2}, \cdots u_{k}$ can then only take on the values $A_{ik}$  in lemma~\ref{abp}.

The conditional expectation $E_{P}( h(\nu) \mid  u_{1}, u_{2}, \cdots u_{k}, B_{t})(\nu)$  is a function defined at the points $\nu_{1}$ and $\nu_{2}$ depending on whether $u_{1}, u_{2}, \cdots u_{k}$ takes on the value of $A_{1k}$ or $A_{2k}$. 
The following result describes this conditional expectation in the present setting.
\begin{lemma}   \label{btx}      \begin{equation}
E_{P}( h(\nu) \mid  u_{1}, u_{2}, \cdots u_{k}, B_{t})=E_{P}(h(\nu) \mid B_{t}) =\frac{1}{2} h(\nu_{1}) + \frac{1}{2} h(\nu_{2}).
\end{equation}
The conditional expectation is therefore constant at the $\nu_{i}$ with value $\frac{1}{2}\, h(\nu_{1}) + \frac{1}{2}\, h(\nu_{2}).$
\end{lemma}
Proof:  In the argument below we denote the vector $u_{1}, u_{2}, \cdots u_{k}$ by $U_{k}$.
For $i=1,2$ and $j=1,2$, apply the independence of lemma~\ref{abp} (for the second equality) to get  the following chain of equalities for conditional probability \begin{equation}
\pi_{i}( h(\nu)= h(\nu_{j}) \mid U_{k}, B_{t})   =  \frac{ \pi_{i}( [h(\nu)= h(\nu_{j})], U_{k}, B_{t}) }{\pi_{i}(U_{k}, B_{t}) } = \end{equation}
\[ \frac{\pi_{i}([h(\nu)=h(\nu_{j})], B_{t}) \cdot \pi_{i}(U_{k} )}{ \pi_{i}( B_{t}) \cdot \pi_{i}(U_{k})} = 
 \pi_{i}(h(\nu)= h(\nu_{j}) \mid B_{t}) \]
 where $\pi_{i}(h(\nu)= h(\nu_{j}) \mid B_{t})$ =1 or 0 according as $i=j$ or $i  \neq  j$.\\
 This implies
\begin{equation}
E_{P}( h(\nu) \mid  U_{k}, B_{t}) =   \frac{1}{2} E_{\pi_{1}}( h(\nu) \mid U_{k}, B_{t})  +  \end{equation}
\[ \frac{1}{2} E_{\pi_{2}}( h(\nu) \mid U_{k}, B_{t} ) =   \frac{1}{2} E_{\pi_{1}}(h(\nu) \mid B_{t})+ \]
\[\frac{1}{2} E_{\pi_{2}}(h(\nu) \mid B_{t}) = \frac{1}{2} h(\nu_{1}) + \frac{1}{2} h(\nu_{2}) = E_{P}( h(\nu) \mid B_{t} ). \]    {\bf Q.E.D.}
    \begin{corollary} $h(\nu_{1}) =h(\nu_{2})$.  \end{corollary}
Proof: The set $\{A_{ik}, B_{t}\}= \{A_{ik} \cap B_{t}\}$ can be described in terms of the given variables and set of the conditional expectation
$E_{P}( h(\nu) \mid u_{1}, u_{2}, \cdots u_{k}, B_{t})$ so by the definition of conditional expectation (see item 7 in the list in section 2.3.2)
\[ \int_{\{A_{ik}, B_{t}\}} E_{P}( h(\nu) \mid u_{1}, u_{2}, \cdots u_{k}, B_{t})\, dP= \int_{\{A_{ik}, B_{t}\}} h(\nu)\, dP. \]
Note that $P(A_{ik}, B_{t})=P(\nu_{i})>0$. Each integrand  in this relation reduces to a constant on $\nu_{i}$; on the left it is  $\frac{1}{2} h(\nu_{1}) +\frac{1}{2} h(\nu_{2})$ by lemma~\ref{btx} and on the right  it is $h(\nu_{i})$. Therefore
\begin{equation}  \frac{1}{2} h(\nu_{1}) +\frac{1}{2} h(\nu_{2}) =  h(\nu_{1})= h(\nu_{2}). \end{equation}
and this concludes the proof of the corollary and therefore of theorem~\ref{indep2}. Then $a=1-a$ in section 2.2 and $a=1/2$. This concludes the proof of theorem~\ref{snb2}.

  \noindent {\bf Acknowledgement} I want to thank  Professor Laurent Moret-Bailly for alerting me to a gap in a proof in a previous version of this work posted on the Math ArXiv and for very helpful criticisms and comments. 

\vspace{.2in}

{\bf Emeritus Professor}  

{\bf Lehman College and Graduate Center, CUNY}

{\bf {\it email:} richard.isaac@lehman.cuny.edu}
\end{document}